\documentclass{article}
\usepackage{diagrams}
\usepackage{amssymb}
\usepackage{amsmath, amscd}
\newcommand{\pf}{{\sc Proof. }}

\newcommand{\LL}{{\mathbb L}}

\newcounter{prop}

\newtheorem{prop}[prop]{Proposition}
\newtheorem{thm}[prop]{Theorem}

\numberwithin{prop}{section}
\numberwithin{subsection}{section}

\title{Projective lines as groupoids with projection structure}

\author{Anders Kock\\
\small Dept.\ of Math., University of Aarhus, Denmark \\
 \normalsize}
\date{}

\begin{document}
\maketitle

\begin{abstract}The coordinate projective line over a field is seen 
as a groupoid with a further `projection' structure. We investigate 
conversely to what extent such an, abstractly given, groupoid may be coordinatized by a 
suitable field constructed out of the geometry.
\end{abstract}

\section*{Introduction}
Given a field $k$ and a 2-dimensional vector space $V$ over it, the 
corresponding projective line $P(V)$ is the set of 1-dimensional 
linear subspaces of $V$. However, it is not a structureless set, but 
rather, it is the set of objects of a groupoid $\LL (V)$, whose morphism are 
the linear isomorphisms between these subspaces. If $A$ and $B$ are 
distinct 1-dimensional subspaces of $V$, a linear isomorphism $A\to 
B$ consists in projection from $A$ to $B$ in a the direction of a 
unique 1-dimensional linear subspace $C$, distinct from $A$ and $B$. 
And vice versa, a 1-dimensional linear subspace $C$ distinct from $A$ 
and $B$ gives rise to such a linear isomorphism.

This provides the groupoid $\LL (V)$ with a certain combinatorial 
structure which we call a {\em projection structure}. The aim of the 
present article is to describe properties of abstract groupoids $\LL$ with such 
structure that will imply that $\LL$ is of the form $\LL (V)$ for 
some vector space $V$ over some field $k$. Note that the field $k$ is 
to be {\em constructed}, unlike in our [Kock, 2010a] where the field  is presupposed. 

The properties of $\LL $  to be described are presented in terms of 
four axioms, whose validity for the case of $\LL (V)$ is easily 
verified algebraically, but to a certain extent have visible geometric 
content, for the case where $k$ is the linear continuum; see e.g.\ 
the illustration to Axiom 4 below. 

\medskip

The present note is a completed form of the preprint [Kock, 2010b], and is a 
continuation of [Kock, 2010a]. The step from  [Kock, 2010a] to the present note 
corresponds to going to the ``deeper level'' described by  [Baer, 1952], Appendix 
I, where he says that the problem of coordinatizing projective 
geometry by means of linear algebra ``{\em can be attacked on two 
essentially different levels: a deeper one where the projective 
geometry is given abstractly, where nothing is supposed to be known 
concerning the underlying linear manifold} [vector space];  {\em and the 
considerably simpler problem where the projective geometry is given 
as the projective geometry of subspaces of some definite linear 
manifold} \ldots.

Baer notes (loc.cit.\ III.4) that ``{\em a line \ldots has no 
geometrical structure if considered  as an isolated or absolute 
phenomenon \ldots. But it receives a definite geometric structure, if 
embedded into a linear manifold of higher rank.} Our contention here 
is that the geometric structure of the (projective) line
 can be described without explicit reference
 to such ``higher rank'' manifold, (the projective plane or 3-space), 
namely by considering it from a (geometrically natural) groupoid theoretic viewpoint. 
This viewpoint was (to our knowledge) first considered in [Kock, 1974], and 
developed further by Diers and Leroy [Diers-Leroy, 1994] and in [Kock, 2010a], 
and, independently, in [Cathelineau, 1995]. 

\section{Projective lines in linear algebra}\label{sect1x}
We collect here some classical facts  from projective geometry when 
seen as an aspect of linear algebra. The emphasis is on the, less 
classical, groupoid theoretic aspect [Kock, 1974], [Diers-Leroy, 
1994], [Cathelineau, 1995], [Kock, 2010a], 
 as alluded to in the Introduction. 

If $A$ and $B$ are mutually distinct 1-dimensional linear subspaces of a 
2-dimensional vector space $V$, and if $\gamma: A \to 
B$ is a $k$-linear isomorphism, the set of vectors in $V$ of the 
form $a - \gamma (a)$ for $a\in A$ form a 1-dimensional subspace $C$ 
distinct from $A$ and $B$. Conversely, a 1-dimensional linear subspace 
$C$  distinct from $A$ and $B$ gives rise to such an isomorphism 
$\gamma$, with $\gamma (a)$ defined as the unique $b\in B$ with 
$a-b\in C$. This $\gamma$, we denote $C:A\to B$ or
$$\begin{diagram}A&\rTo^{C}&B\end{diagram}.$$
This establishes a bijective correspondence between linear 
isomorphisms $A\to B$, and 1-dimensional linear subspaces $C$ distinct 
from $A$ and $B$.   

We present some properties of this bijective correspondence, and the 
groupoid structure; we do it in a fully coordinatized situation, so we 
assume that the vector space $V$ is $k^{2}$
where $k$ is a field. It gives rise to a model $\LL (k^{2})$  of the axiomatics 
presented, namely the projective line over $k$.  Its set of objects 
is the set of points  
$P(k^{2})$. A point $A$ in $P(k^{2})$ is a 1-dimensional 
linear subspace  $A \subseteq k^{2}$, and it may be presented by any 
non-zero vector $a=(a_{1},a_{2})\in A$ with $a_{1}$ and $a_{2}$ $\in 
k$ (in traditional notation ``$A$ is the point with homogeneous 
coordinates 
$(a_{1}:a_{2})$''). So $a$  is a basis vector  for $A$.
The calculation of $C:A\to B$ in coordinates is the following. Here, $A,B$, and $C$ are 
1-dimensional subspaces of the 2-dimensional vector space $k^{2}$. Pick  basis vectors $a$ 
for $A$, $b$ for $B$ and $c$ for $C$.  Then the linear map $C:A\to B$ 
has a matrix w.r.to the bases $\{a\}$ for $A$ and $\{b\}$ for $B$; it 
is a $1\times 1$ matrix, whose only entry is the scalar
$$\frac{|a \; c|}{|b \; c]}$$
where $|a \; c|$ denotes the determinant of the $2\times 2$ matrix 
whose columns are $a\in k^{2}$, $c\in k^{2}$, and similarly for $|b\; 
c|$. Since composition of 
linear maps corresponds to multiplication of matrices, which here is 
just multiplication of scalars, it follows that  the composite map
\begin{equation}\label{birapx}\begin{diagram}A&\rTo ^{C}&B&\rTo^{D}&A
\end{diagram}\end{equation}
is (multiplication by) the scalar
\begin{equation}\label{birapx2}\frac{|a \; c|}{|b \; c]}\cdot \frac{|b \; d|}{|a \; 
d]}\end{equation}
where $d$ is any basis vector for $D$; and this scalar is independent 
of the choice of the four basis vectors $a,b,c,d$ chosen, and it is 
the classical expression for cross ratios in $P(k^{2})$. (This was 
also argued geometrically in [Kock, 1974]  p.\ 3, and algebraically 
in [Diers-Leroy, 1994], 
Theorem 1-5-3.)  The cross ratio (\ref{birapx})  is  
denoted $(A,B;C,D)$. 
-- Note that the expression (\ref{birapx2}) does not depend on $c$ and $d$ being 
linearly independent.

We note the following two trivial equations, which will be referred 
to later:
\begin{equation}\label{ax1x}\frac{|a \; c|}{|b \; c]}\cdot \frac{|b 
\; c|}{|a \; 
c]}=1; \quad \quad \frac{|a \; c|}{|b \; c]}\cdot \frac{|b \; c|}{|d \; 
c]}=\frac{|a \; c|}{|d \; c]}.
\end{equation}
This does not use any properties of determinants. Using that 
determinants are alternating, $|b \; c| = -|c\; b|$, we have also
the equality
\begin{equation}\label{ax2x}\frac{|a \; c|}{|b \; c|}\cdot \frac{|b \; d|}{|a\; d]}\cdot 
\frac{|a \; b|}{|c \; b]}=\frac{|a \; b|}{|c \; b]}\cdot \frac{|c \; 
a|}{|d \; a]}\cdot 
\frac{|d \; b|}{|c \; b]},\end{equation}
 because the factors on the left may be paired off 
with factors on the right 
(modulo four sign changes, which cancel). Similarly, the factors
  on the left in the following equation may be paired off modulo 
three sign changes, 
whence the 
minus sign:
\begin{equation}\label{ax4x}\frac{|a \; c|}{|b \; c|}\cdot \frac{|b \; 
a|}{|c\; a]}\cdot 
\frac{|c \; b|}{|a \; b]} = -1\end{equation}
Finally, we have the following well known relationship between 
cross-ratios 
\begin{equation}\label{ax3x} 1- (A,B;C,D) = (A,C;B,D)
\end{equation}
in coordinates (with $a =(a_{1},a_{2})$ a basis vector for $A$, 
and similarly for $b,c,d$) the two sides to be compared are 
$$1- (A,B;C,D)=  1-  \frac{|a \; c|}{|b \; c]}\cdot \frac{|b \; d|}{|a \; d]} 
= \frac{|b\; c|.|a\; d| - |a\; c|.|b\; d|}{|b\; c|.|a\; 
d|},$$
and $$(A,C;B,D) = \frac{|a\; b|}{|c\; b|}\cdot \frac{|c\; d|}{|a\;d|} 
= \frac{|a\; b|.|c\; d|}{|c\; b|.|a\;d|}.$$
The denominators in these two expressions are the same except for 
sign, so it suffices to see that the enumerators are the same except 
for sign. Using $|a\; b| = a_{1}b_{2}-a_{2}b_{1}$ etc., we get for 
the first enumerator
$$(b_{1}c_{2}-b_{2}c_{1}).(a_{1}d_{2} - a_{2}d_{1})- 
(a_{1}c_{2}-a_{2}c_{1}).(b_{1}d_{2}-b_{2}d_{1}),$$
and multiplying out, we get eight terms; there are two terms 
$a_{1}b_{1}c_{2}d_{2}$ with opposite sign, and similarly two terms 
$a_{2}b_{2}c_{1}d_{1}$, and these terms cancel, The four 
remaining terms are then seen to be the terms of the enumerator in 
the expression for $(A,C;B,D)$, except for sign. This proves 
$1-(A,B;C,D) = (A,C;B,D)$, as desired.

\section{The axiomatics}\label{threefourx}
A {\em projection structure} on a groupoid $\LL$ consists in the 
following: for any two distinct objects $A$ and $B$ of $\LL$, 
there is given a bijection between the set $\hom_{\LL}(A,B)$ and the 
set $|\LL|\backslash \{ A,B\}$, where $|\LL|$ denotes the set of 
objects of $\LL$. The arrow $A\to B$ corresponding to an object $C$ 
(with $C\neq A,B$) is denoted $C:A\to B$ or 
$$\begin{diagram}A&\rTo^{C}&B.
\end{diagram}$$
A functor $P: \LL \to \LL'$ between groupoids with projection 
structure is said to {\em preserve} these structures if it is 
injective on objects, and if 
$$\begin{diagram}[midshaft]P(A&\rTo^{C}&B ) \quad = \quad  P(A) 
&\rTo^{P(C)}& P(B).
\end{diagram}$$
Note that the injectivitity condition is needed to make sense to the 
right hand expression.

A functor $P$ preserving projection structures in this sense is also 
called a {\em homomorphism} or a {\em projectivity} (whence the letter 
$P$). It is easy to see 
that a homomorphism, as described, is a {\em faithful} 
functor $\LL \to \LL'$.

If $\LL$ with a projection structure has only one or only two objects, 
it is a group, respectively a disjoint 
union of two groups, and in this case, a projection structure is a void concept. So 
henceforth, we assume that $\LL$ has at least three distinct objects.

If a groupoid $\LL$ with a projection structure has at least three objects, 
it is {\em connected} in the sense that for any pair $A,B$ of 
objects, $\hom _{\LL}(A,B)$ is nonempty: any object $C \neq A,B$ 
defines a morphism $C:A\to B$.

In a connected groupoid $\LL$, all vertex groups $\hom _{\LL}(A,A)$ 
are isomorphic: conjugation by an arrow $A\to B$ provides an 
isomorphism $\hom _{\LL}(A,A)\to \hom _{\LL}(B,B)$. If one of these 
vertex groups is commutative, all vertex groups are commutative, and 
the isomorphism  $\hom _{\LL}(A,A)\to \hom _{\LL}(B,B)$ defined by 
conjugation by an arrow $\alpha :A \to B$ does not depend on the 
choice of the arrow $\alpha :A \to B$. For brevity, we call such (connected) groupoids 
{\em commutative}. In such a groupoid, all vertex groups $\LL 
(A,A)$ are {\em canonically} isomorphic, and may be identified with 
one single (commutative) group, the {\em group of abstract scalars} 
of $\LL$; the group $\hom _{\LL}(A,A)$ may be called {\em the group 
of scalars at $A$}; it is canonically isomorphic to the group of 
abstract scalars. If $\mu \in \hom_{\LL}(A,A)$ and $\mu' \in 
\hom_{\LL}(B,B)$ represent the same abstract scalar, we write $\mu 
\equiv \mu'$; this is the case if there is an arrow $\alpha :A\to B$ 
conjugating $\mu$ to $\mu'$, equivalently, with $\mu.\alpha = \alpha 
.\mu'$. (We compose from left to right in $\LL$.) Sometimes we do not 
distinguish between abstract scalars, and scalars at $A$, and between 
$=$ and $\equiv$ for scalars. 

\medskip

We present some axioms for groupoids $\LL$ with projection structure, 
with $\LL$ being a commutative groupoid with at least three objects.
The objects of $\LL$, we shall also call {\em points} of $\LL$, 
because of the intended geometric interpretation of $\LL$.

\medskip

\noindent{\bf Axiom 1.}
{\em The following two diagrams commute, where
$A,B,C$, and $D$ are  mutually distinct points of $\LL$: 
$$\begin{diagram}[nohug]
A&\rTo^{C}&B\\
&\rdTo_{1_{A}}&\dTo_{C}\\
&&A\\
\end{diagram}
\quad \quad \mbox{ and } \quad 
\begin{diagram}[nohug]A&\rTo^{C}&B\\
&\rdTo_{C}&\dTo_{C}\\
&&D.
\end{diagram}.$$}

The validity of this axiom in the coordinate model $P(k^{2})$ is 
expressed by (\ref{ax1x}).

\medskip

We shall introduce the 
notion of {\em cross ratio}, and a matrix 
notation for it. Cross ratios are certain 
scalars. 
Consider  four distinct points  
$A,B,C,D$ of $\LL$.
We write $$(A,B;C,D) \mbox{ \quad or \quad }  \left[ 
    \begin{array}{cc}A&B\\
	C&D
    \end{array}\right]$$
    for the scalar at $A$, or for the abstract scalar it represents, given 
    as the composite
    $$\begin{diagram} A&\rTo ^{C}& B&\rTo ^{D}& A.
    \end{diagram}$$
and call it the {\em cross ratio} of the 4-tuple. The word {\em 
bi-rapport} will also be used.

It is an immediate consequence of the notion of projectivity, i.e.\ a morphism $P$ of groupoids 
with projection structure,  that such projectivity $P$ preserves cross ratios
$$(P(A),P(B);P(C),P(D))= P((A,B;C,D)).$$ 

We note that if 
$\LL$ has only three points, there are no scalars that appear as 
cross ratios. On the other hand, if there are more than three points, 
every scalar except $1$ appears as a cross ratio, more precisely
\begin{prop}\label{crrxx}Given $\mu \in \LL (A,A)$ with $\mu \neq 1_{A}$.
 For any $B$ and $C$, mutually 
distinct, and distinct from $A$, there is a unique $D$, distinct from 
$A$, $B$,  and $C$, such that $\mu = (A,B;C,D)$.
\end{prop}
\pf Since $\LL$ is a groupoid, 
there is a unique  arrow $d$ making the triangle 
$$\begin{diagram}[nohug]A&\rTo ^{C}&B\\
&\rdTo_{\mu}&\dTo_{d}\\
&&A
\end{diagram}$$
commute. Since $B$ and $A$ are distinct, $d$ is of the form $D:B\to 
A$ for a unique point $D$ distinct from $A$ and $B$; and, by Axiom 1, 
$C=D$ would  imply that we get  $\mu = 1_{A}$, which we excluded.

\medskip

In a cross ratio 
expression $(A,B;C,D)$, one might allow the possibility that $C=D$, 
in which case the first part of Axiom 1 would say that $(A,B;C,C)= 
1_{A}$; however, this usage would make some of the statements to be made 
later more clumsy. This is why we consider only cross ratio 
expressions with four distinct entries.

\medskip

\noindent{\bf Axiom 2.} {\em  The following  diagram commutes, where
$A,B,C$, and $D$ are  mutually distinct points of $\LL$: 
$$\begin{diagram}[nohug]
A&\rTo ^{C}&B&\rTo ^{D}&A\\
    \dTo ^{B}&&&&\dTo _{B}\\
    C&\rTo _{A}&D&\rTo _{B}&C.
    \end{diagram}$$}

The validity of this axiom in the coordinate model $P(k^{2})$ is 
expressed by (\ref{ax2x}).

\begin{prop}\label{cyclicx} In a cross ratio matrix, the columns may be 
interchanged, without changing the absolute scalar which the matrices 
represent. Also, the rows may be interchanged without 
changing the absolute scalar which the matrices 
represent.
\end{prop}
\pf   For columns:
Consider the diagram
    $$
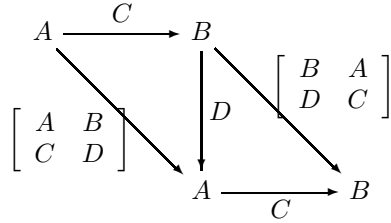
\begin{diagram}[nohug]A&\rTo ^{C}&B&&\\
    &\rdTo _{\left[ \begin{array}{cc}A&B\\C&D
    \end{array}\right] }&\dTo _{D}&
    \rdTo ^{\left[ \begin{array}{cc}B&A\\D&C
    \end{array}\right] }&\\
    &&A&\rTo_{C}&B
    \end{diagram}$$
    The two triangles commute by definition; hence, the  total quadrangle
commutes; but this can be expressed:  $C:A\to B$ conjugates the scalar
    $(A,B;C,D)$ at $A$ to the scalar $(B,A;D,C)$ at $B$. For rows, 
this is just a reformulation of Axiom 2, which in terms of cross 
ratios says that $B:A\to C$ conjugates $(A,B;C,D)$ to $(C,D;A,B)$.

\begin{prop}\label{assxx}We have
$$\left[ \begin{array}{cc}
A&B\\ D&C 
\end{array} \right] = \left[ \begin{array}{cc}
A&B\\
C&D
\end{array} \right]^{-1}$$
and
$$\left[ \begin{array}{cc}A&B\\ C&D
\end{array} \right].\left[ \begin{array}{cc}
A&B\\ D&E
\end{array} \right] = \left[ \begin{array}{cc}A&B \\ C&E
\end{array} \right].$$
\end{prop}
\pf  For the first assertion, we consider the composite
$$\left[ \begin{array}{cc}
A&B\\ C&D 
\end{array} \right] .\left[ \begin{array}{cc}
A&B\\
D&C
\end{array} \right]$$
which amounts to a four-fold composite
$$\begin{diagram}A&\rTo^{C}&B&\rTo^{D}&A&\rTo^{D}&B&\rTo^{C}&A
\end{diagram}.$$
The two middle factors give $1_{B}$, by Axiom 1, and then the two 
remaining factors give $1_{A}$, again by Axiom 1.  -- The proof of 
the second assertion is similar, just replace the last occurrence of 
$C$ by $E$; then the two middle factors again give $1_{B}$, and then 
the two remaining factors give the defining composite for$(A,B;C,E)$.

\medskip

We now state the third axiom which refers to ``middle four 
interchange'' in cross ratio expressions like $(A,B;C,D)$:

\medskip

\noindent{\bf Axiom 3.} {\em If $(A,B;C,D)\equiv (A',B';C',D')$, then also
$(A,C;B,D)\equiv (A',C';B',D')$.}

\medskip

\begin{sloppypar}The validity of this axiom in the coordinate model 
$P(k^{2})$ follows from
 
(\ref{ax3x}), since if $(A,B;C,D)=(A',B';C',D')$, $=\mu$, say, then
$(A,C;B,D)= 1-\mu = (A',C';B',D')$,  by two applications of 
(\ref{ax3x})).\end{sloppypar}
 
\medskip

Middle four interchange is clearly an involution.
Since every scalar $\mu \neq 1$ appears as a cross ratio of four
distinct points, it follows from the Axiom that we have:
\begin{prop}\label{mfx} There is an 
involution $\Phi$ on the set $G\backslash \{1\}$ of scalars $\neq 1$, 
such that for any distinct 4-tuple $(A,B,C,D)$, we have
$$(A,C;B,D) = \Phi (A,B;C,D).$$
\end{prop}

\medskip
 
In a projective line over a field, $\Phi (\mu ) = 1-\mu$, and this 
involves the additive structure (in terms of the binary operation ``minus'') of the field, 
which is not assumed in our context, but is rather something to be 
constructed.

\medskip

We note that a projectivity  $P: \LL \to \LL'$ of groupoids with 
projection structure, both of which  satisfy Axioms 1-3, preserves the 
corresponding involutions $\Phi$ and $\Phi'$: present a given scalar 
$\mu \neq 1$ of $\LL$ as $(A,B;C,D)$; then $\Phi (\mu ) = (A,C;B,D)$, 
and so
$$P(\Phi (\mu ))=P(A,C;B,D) = (P(A),P(C);P(B),P(D)) $$ $$= \Phi' 
(P(A),P(B);P(C),P(D))= 
\Phi' (P(\mu )).$$
Thus the group of scalars is a $\Phi$-group, in the following sense:

\medskip

\begin{sloppypar}By a {\em $\Phi$-group}, we understand a (multiplicatively written) 
commutative group $G$ equipped with an involution $\Phi$ on the set 
$G\backslash \{1\}$; a morphism of such is an injective group 
homomorphism compatible with the $\Phi$s.
\end{sloppypar}

\medskip

\noindent{\bf 2.5 \; Permutation laws. }\stepcounter{prop}
Let $c:Q\to H$ be a surjection. If $\sigma :Q\to Q$ is an endo-map 
(we are only interested in permutations),
then we say that $\sigma$ {\em descends} along $c$ if there exists an 
endomap $\overline{\sigma}: H\to H$ (necessarily unique, by 
surjectivity of $c$) such that the square
$$\begin{diagram}Q&\rTo^{c}&H\\
\dTo^{\sigma}&&\dTo_{\overline{\sigma}}\\
Q&\rTo^{c}&H
\end{diagram}$$
commutes.  A composite of two endomaps which descend along 
$c$ again descends along $c$. Clearly, $\sigma$ descends along $c$ iff
$c(q)=c(q')$ implies $c(\sigma (q)) = c(\sigma (q'))$ for all $q$ and 
$q'$ in $Q$.

We have in mind the case where $Q$ is the set $|\LL |^{<4>}$ of 4-tuples of mutually 
distinct points of $\LL$, and where $c$ is cross ratio formation, 
with $H$ the set of (abstract) scalars  $\neq 1$ of $\LL$.
Then the following proposition distills some information already 
obtained or assumed; $\LL$ is a groupoid with projection structure, 
satisfying Axioms 1-3.

\begin{prop}\label{permxx}Every permutation $\sigma$ of  $|\LL |^{<4>}$ 
arising from a permutation $\sigma \in \mathfrak{S}_{4}$ descends 
along the cross ratio formation map $c: |\LL |^{<4>}\to G\backslash 
\{1\}$.
\end{prop}
\pf  The 24 permutations on $|\LL |^{<4>}$ (whose 
elements we, as above, write as $2\times 2$ matrices) arising from 
$\mathfrak{S}_{4}$ are generated by the following four permutations:
interchange of rows, intechange of columns, interchange of two lower 
entries, and middle four interchange. Interchange of rows, and 
interchange of columns descend to the identity map of $G\backslash 
\{1\}$, by Proposition 
\ref{cyclicx}; 
interchange of two lower entries descends 
to multiplicative inversion $(-)^{-1}$, 
by (the first assertion of) Proposition \ref{assxx}; and, finally, middle four 
interchange descends to $\Phi$, by Proposition \ref{mfx}.

\medskip

In a similar vein, let $p:G \to G'$ be a $\Phi$-group homomorphism, let $P: 
|\LL | \to |\LL'|$, and let let $p:G \to G'$ be a $\Phi$-group 
homomorphism, where $G$ and $G'$ are the groups of scalars of $\LL$ 
and $\LL'$, respectively.
 We assume  Axioms 1-3.  Then
\begin{prop}\label{pex} If $P$ preserves the cross ratio of 
a given distinct 4-tuple (relative to $p$), then it also preserves any permutation 
instance of it.
\end{prop}
\pf  It suffices to prove it for each of the four generating 
permutations mentioned above. The four cases are proved along the 
same lines, so let us just give the proof for the case of middle four 
interchange. So assume that a cross ratio $(A,B;C,D)$ is preserved, so
 $p((A,B;C,D))= (P(A),P(B);P(C),P(D))$.
Then 
$$p((A,C;B,D)) = p(\Phi (A,B;C,D)) = \Phi '(p(A,B;C,D))=$$ $$
=\Phi'(P(A),P(B);P(C),P(D))
= (P(A),P(C);P(B),P(D))$$
where $\Phi$ and $\Phi'$ denote the involution arising from middle 
four interchange in $\LL$ and $\LL'$, respectively.

\medskip

\noindent{\bf 2.8 \;  Tri-rapports, and the scalar $-1$} 
\stepcounter{prop}
Cross ratios in projective geometry are also  
called {\em bi-rapports}. In our context, this is quite logical, 
since a cross ratio is an endomap presented as a composite of two 
maps (which themselves are not endomaps). Similarly, a tri-rapport in 
our context is an endomap presented as a composite of three 
(non-endo) maps,
$$\begin{diagram}
A&\rTo^{X}&B&\rTo^{Y}&C&\rTo^{Z} &A
\end{diagram},$$ and this scalar (at $A$) may be denoted in matrix 
form, in analogy with cross ratios (bi-rapports)
$$\left[ \begin{array}{ccc}
A&B&C\\
X&Y&Z \end{array}\right].$$
 This usage (in our context) is found is e.g.\ [Diers-Leroy, 1994] (even 
generalized  into $n$-rapports), or (vastly 
generalized) in [Cathelineau, 1995]. We may also 
denote tri-rapports in ``one-line'' form, by $(A,B,C;X,Y,Z)$, to save space.
\begin{prop}\label{cycltrix}The columns in a tri-rapport matrix may be permuted 
cyclically without changing its value.
\end{prop}
\noindent The proof is analogous to the proof of Proposition \ref{cyclicx} (the 
first part), and is left to the reader.

\medskip

Of particular importance are tri-rapports of the following form
\begin{equation}\label{trirx}\begin{diagram}
A&\rTo^{C}&B&\rTo^{A}&C&\rTo^{B} &A
\end{diagram}\end{equation}
so, in the matrix notation, the lower row is a cyclic permutation of 
the upper row. The geometric content of this scalar, in the classical 
geometric model of the real projective line, is that it is the scalar 
$-1$. For, consider the picture (from [Kock, 1974])

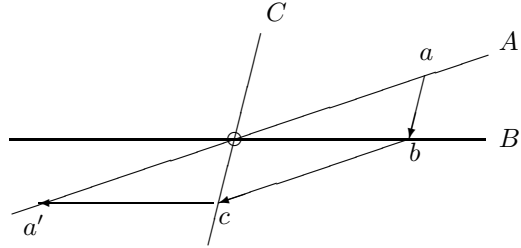
\begin{picture}(50,110)(-195,-50)
    \put(-10,-40){\line(1,4){20}}
    \put(0,0){\circle{5}}
    \put(-85,0){\line (1,0){180}}
    \put(-84,-28){\line(3,1){180}}
    \put(12,45){$C$}
    \put(100,34){$A$}
    \put(100,-3){$B$}
    \put(72,24){\vector(-1,-4){6}}
    \put(70,30){$a$}
    \put(66,0){\vector(-3,-1){72}}
    \put(66,-8){$b$}
    \put(-8,-24){\vector(-1,0){66}}
    \put(-6,-32){$c$}
    \put(-80,-36){$a'$}

\end{picture}

\noindent Then given an element $a\in A$, applying consecutively the linear maps $C:A\to 
B$, $A:B\to C$, and $B:C\to A$ 
to $a$, as illustrated, gives reflection of $a$ in $O$. But 
such reflection is is independent 
of the choice of $B$ and $C$. In our context, the reflection is an 
arrow $A\to A$, that is, a scalar at $A$, and we would like to 
denote it $(-1)_{A}$; to do this we need that it is independent of 
the choice of $B$ and $C$. 
This independence seems not to be something that we can derive on the basis of the axioms 
given so far, so we need this as a further axiom. We shall see that it suffices to 
require independence of the choice of $C$, so we pose:

\medskip

\noindent {\bf Axiom 4.} {\em Given four distinct points $A,B,C,D$,
 we have
$$\left[\begin{array}{ccc}A&B&C\\
C&A&B \end{array} \right] = 
\left[\begin{array}{ccc}A&B&D\\ D&A&B \end{array} \right].$$}

The validity of this axiom in the coordinate model $P(k^{2})$ follows from 
(\ref{ax4x}), since both sides of the equation in this case give 
$-1$.

\medskip 

We can then deduce
\begin{prop}\label{eightx} The scalar at $A$ given by $(A,B,C;C,A,B)$
 is independent of the choice of $B$ and 
$C$, and may thus be denoted $(-1)_{A}$. Furthermore, $(-1)_{A} 
\equiv (-1)_{B}$ for any $A$ and $B$, so they represent the same 
abstract scalar, denoted $(-1)$. Finally, $(-1).(-1) = 1$.
\end {prop}
\pf  The independence of the choice of $C$ is a reformulation 
of the Axiom. The independence of the choice of $B$ then follows purely 
formally from Proposition \ref{cycltrix}. To prove $(-1)_{A}\equiv 
(-1)_{B}$, pick  $C$ distinct from $A$ and $B$.  By what is already 
proved, $(-1)_{A}$ may be presented by $(C:A\to B).(A:B\to C).(B:C\to 
A)$, and $(-1)_{B}$ may be presented by 
$(A:B\to C).(B:C\to A).(C:A\to B)$. Now consider the four-fold composite
$$\begin{diagram}A&\rTo^{C}&B&\rTo^{A}&C&\rTo^{B}&A&\rTo^{C}&B
\end{diagram};$$ the composite of the three first of these maps is 
$(-1)_{A}$ and the composite of the three last is $(-1)_{B}$. Then 
the two ways of interpreting the four-fold composite, by the 
associative law, expresses that $C:A\to B$ conjugates $(-1)_{A}$ to 
$(-1)_{B}$, i.e.\  expresses $(-1)_{A} \equiv (-1)_{B}$. The fact 
that $(-1).(-1)=1$ follows if we can prove that $(-1)_{A}.(-1)_{A }$ 
is the identity map at $A$. This follows by three apllications of 
Axiom 1, by representing the first factor $(-1)_{A}$ by the tri 
rapport
$(A,B,C;C,A,B)$, and the second one by the tri-rapport 
$(A,C,B;B,A,C)$.

\medskip

	\noindent{\bf Remark.} Assume $-1 \neq 1$. Then the classical way of dealing with the scalar $-1$, here 
	defined as a tri-rapport, is in terms of {\em harmonic 
	conjugates}:	         
  Given $A,B,C, H$. Then 
$$(A,B;C,H)= -1 $$ iff $H:B\to A$ equals the composite
$$\begin{diagram}
B&\rTo^{A}&C&\rTo^{B}&A\end{diagram}.$$
For, precomposing the composite with $C:A\to B$ gives $-1$, and 
 precomposing $H:B\to A$ with $C:A\to B$ gives $-1$ iff  $(A,B;C,H)=-1$.
 The classical  way of formulating this 
characterizing property of $H$  is: $H$ is the 
{\em harmonic conjugate} of $C$ w.r.to $A,B$.

\medskip

 For $\mu$ a scalar, $-\mu$ denotes 
of course $(-1).\mu$ ($ = \mu . (-1)$). We note the following sign change property of certain tri-rapports 
(with only four points occurring);.
\begin{prop}\label{signchangex} Let $A,B,C,D$ be mutually distinct points. Then
$$(A,B,C;D,A,D) = - (A,B,D;D,A,C)$$
\end{prop}
\pf 
Consider the diagram
$$\begin{diagram}
A&\rTo^{D}&B&\rTo^{A}&D&\rTo ^{C}&A\\
\dTo^{1}&&&&\dTo_{A}&&\dTo_{-1}\\
A&\rTo_{D}&B&\rTo_{A}&C&\rTo_{D}&A
\end{diagram}.$$
The left hand rectangle commutes (insert a ``North-East'' arrow $A:B\to D$ and use 
Axiom 1); the right hand square commutes by (a variant of) the 
definition of the scalar $-1$.

\medskip

\noindent{\bf Remark.} One cannot conclude that $-1$ is distinct from 
$1$; take $\LL (k^{2})$ for $k$ a field of characteristic 2. 
I don't know whether a four-point $\LL$ has $-1 \neq 1$.
The projective 
line $P(k^{2})$ over the field $k$ with three elements has $-1\neq 
1$, and is up to isomorphism the only such groupoid with projection 
structure, satisfying the axioms 1-4 and $-1\neq 1$. For, it can be seen by a 
straightforward sudoku argument that such $\LL$ is unique up to 
isomorphism: If the four points are $A,B,C,D$, 
the composition is uniquely determined; for instance
$(C:A\to B).(A:B\to C)$ is necessarily $D:A\to B$; for, it cannot be 
$B:A\to C$, since post-composing with $B:C\to A$ would on the one 
hand give $1$, by Axiom 1, and on the other hand, it would be a 
three-fold composite defining the scalar $-1$. See also the Example 
at the end of [Kock, 2010a]. 

\section{Scalars as 2-cells} 
The purpose of the present Section is to give a proof of  
a relationship  between bi-rapports (cross ratios) and tri-rapports; it 
is (to my knowledge)  first proved in [Diers-Leroy, 1994], Propositon 2-4-3. We present an alternative 
and somewhat more conceptual proof.

 Given a connected commutative 
groupoid $\LL$, with group $G$ of abstract scalars, there is a 
2-dimensional groupoid with $\LL$ as underlying 1-dimensional 
groupoid: namely a 2-cell 
$$\begin{diagram}A&\pile{\rTo^{f}\\ \Downarrow\alpha \\ \rTo_{g}}&B
\end{diagram}$$
is an abstract scalar $\alpha$ such that the scalar $ \alpha_{A}$ at $A$ which 
represents $\alpha$ satisfies $\alpha_{A} .g  = f$, or equivalently, such 
that the scalar $\alpha_{B}$ at $B$ which represents $\alpha$ 
satisfies $g.\alpha _{B} =f$. Both horizontal and vertical composition 
of 2-cells come about from the multiplication in $G$.

For two parallel arrows $f$ and $g$ from $A$ to $B$, there is 
therefore a unique 
2-cell $f\Rightarrow g$, (represented by the scalar $f.g^{-1}$ at $A$ or 
by the scalar  $g^{-1}.f$ at $B$). If $f:A \to A$ is a scalar at $A$, 
the unique 2-cell $ f\Rightarrow 1_{A}$ is the abstract scalar 
associated to $f$.

There is thus nothing in the 2-dimensional structure which is not 
present in the information contained in the 1-dimensional structure, 
so the reader who wants a purely 1-dimensional proof of Proposition 
\ref{bitrix} below is 
referred to [Diers-Leroy, 1994]; however, the result and the proof seem more 
conceptual and compelling in the 2-dimensional formulation.

For the case of a groupoid $\LL$ with projection structure, 
satisfying Axiom 1, we see that the abstract scalar $\mu$ defined by 
the cross ratio $(A,B;C,D)$ (as an arrow $A\to A$) is a 2-cell $\mu$
\begin{equation}\label{cxxx}\begin{diagram}A&\pile{\rTo^{(A,B;C,D)} \\ \Downarrow \mu 
\\ \rTo_{1_{A}}}&A
\end{diagram}\end{equation}
but as an abstract scalar,  
$\mu$ is also 
a 2-cell between two arrows which are not endo-arrows:
$$\begin{diagram}A&\pile{\rTo^{C}\\ \Downarrow\mu \\ \rTo_{D}}&B
\end{diagram}$$
as can be seen by postcomposing $\mu$ in (\ref{cxxx}) by $D:A\to B$ and using Axiom 
1. Identifying scalars at $A$ with abstract scalars, we thus have the 
useful alternative diagrammatic way of defining cross ratios.

\medskip

Consider now a (2-dimensional) diagram
$$\begin{diagram}A&\pile{\rTo^{E}\\ \Downarrow  \ \epsilon \\ 
\rTo_{E'}}&B&\pile{\rTo^{F}\\ \Downarrow \ \phi \\ 
\rTo_{F'}}&C&\pile{\rTo^{G}\\ \Downarrow \gamma \\ 
\rTo_{G'}}&A 
\end{diagram}.$$
The  2-cell obtained by horizontal composition is then the product in 
$G$ of  three cross ratios $\epsilon = (A,B;E,E')$, $\phi =(B,C; F,F')$, and 
$\gamma =(C,A;G,G')$. On the other hand, if the bottom composite is $1_{A}$, the 
(unique) 2-cell  from the top composite to the bottom one is the 
scalar represented by the top composite, i.e.\ the tri-rapport 
$(A,B,C;E,F,G)$. We have therefore proved the following (cf.\  
[Diers-Leroy, 1994],  Prop. 2-4-3):
\begin{prop}\label{bitrix} If $(A,B,C;E',F',G') = 1$, then
$$(A,B,C;E,F,G)= (A,B;E,E').(B,C;F,F').(C,A;G,G').$$
\end{prop}

We note that given distinct $A,B,C$, it is always possible to find 
points $E',F',G'$ such that $(A,B,C;E',F',G') = 1$, provided there 
are at least five points in $\LL$; for, pick distinct $E'$ and $F'$, 
both of them distinct from $A,B$, and $C$, then $(E':A\to B).(F':B\to C)$, as 
a non-scalar, is labelled by a unique $G'$ distinct from $A,C$ and 
from $E',F'$, and then $E',F',G'$ will do the job (postcompose with 
$G':C\to A$ and use Axiom 1).

\begin{prop}\label{bitrixx}Assume that $\LL$ has at least five points, and that
$P: |\LL | \to |\LL '|$, and $p:G\to G'$ are given, such that
$P$ preserves cross ratio formation (w.r.to $p$), with $p$ an 
injective group homomorphism. Then $P$ preserves tri-rapports. 
\end{prop}
\pf  Consider a tri-rapport, say
$(A,B,C;E,F,G)$. Pick, as argued above, $E',F',G'$ with 
$$(A,B,C;E',F',G') = 1.$$ 
 Then applying Proposition \ref{bitrix} for $\LL$, one rewrites the scalar 
$(A,B,C;E,F,G)$ as a pro\-duct of three scalars, each of which is a 
bi-rapport, and hence is preserved by $(P,p)$, i.e.\
$p(A,B;E,E')= (P(A),P(B); P(E),P(E'))$ etc., and since $p$ preserves 
products of scalars, it follows that\newline 
$p(A,B,C;E,F,G)=(P(A),P(B),P(C);P(E),P(F),P(G))$, using now 
Proposition \ref{bitrix} for $\LL'$.

\medskip


\section{Three-transitivity (``Fundamental Theorem'')}
The following Section presents a version of the classical 
``three-transitivity'' theorem for projecti\-vities (called The 
Fundamental Theorem of Projective Geometry in e.g.\ [Coxeter, 1942] 2.8.5). We use Axioms 1-3, 
but not Axiom 4.

Recall that a morphism $P:\LL \to \LL'$ of groupoids with projection 
structure in particular is a faithful functor. Therefore, it gives 
rise to an injective group homomorphism $p:G\to G'$ on the 
corresponding groups of scalars, in particular, it restricts to a map
$p: G\backslash\{1\} \to G'\backslash\{1\}$   
\begin{prop}If $\LL$ and $\LL'$ satsify the 
three axioms 1-3, this map $p$ is compatible with the 
involutions $\Phi$ which we get by virtue of Axiom 3.
\end{prop}
\pf  Given $\mu \neq 1$ in $G$. We have to prove that 
$\Phi'(p(\mu )) = p (\Phi 
(\mu ))$, where $\Phi'$ is the involution on $G'\backslash\{1\}$. 
Since $\mu \neq 1$, we may assume that $\mu = (A,B;C,D)$ for four 
distinct points in $\LL$, and so $\Phi (\mu ) = (A,C;B,D)$, which is 
by definition of cross ratios is the composite
$$\begin{diagram}A&\rTo^{B}&C&\rTo^{D}&A
\end{diagram}$$
in the groupoid $\LL$. Applying the functor $P$ thus gives 
$$p(\Phi (\mu ) ) = 
\begin{diagram}P(A&\rTo^{B}&C)\end{diagram}.\begin{diagram}[midshaft]&P(C\rTo^D &A),
\end{diagram}$$
and since $P$ is assumed to be a morphism of projection structure, this
is the composite in $\LL'$
$$\begin{diagram}P(A)&\rTo^{P(B)}&P(C)&\rTo^{P(D)}&P(A)
\end{diagram}$$
which in turn equals 
$(P(A),P(C);P(B),P(D)) = \Phi'(P(A),P(B);P(C);P(D))$.\\
Now $(P(A),P(B);P(C),P(D))$ is $p(\mu )$, again because $P$ is a 
functor, so that we get $\Phi'(p(\mu 
))$.

\medskip

Let $\LL$ and $\LL'$ be groupoids with projection 
structure, satisfying axioms 1, 2, and 3, and let $G$ and $G'$ be the 
corresponding $\Phi$-groups of scalars. Assume that $\LL$ has at 
least five points.
\begin{thm}\label{main3x}Let $p:G\to G'$ be a (injective) homomorphism of 
$\Phi$-groups. Then 
if $A,B,C$ are distinct points in $\LL$, and  $A',B',C'$ 
are distinct 
points in $\LL'$, there exists a unique morphism of groupoids 
with projection structure $P:\LL \to \LL'$ which on scalars restricts 
to the given $p$ and which satisfies
$P(A)=A', P(B)=B', P(C)=C'$. If $p: G\to G'$ is an isomorphism, then 
so is the morphism $P$. 
\end{thm}
\pf  Part of the proof is from [Kock, 2010a] Section 3, but 
now it is in a 
more general setting. 
The construction, and the uniqueness, of $P$ is forced: 
if $X\neq A,B,C$ in $\LL$, then if $P$ is a morphism, it preserves the cross ratio of 
($A,B;C,X)$, so we are forced to define $P(X)$ to be the unique point 
of $\LL'$ (cf.\ Proposition \ref{crrxx} for $\LL'$) satisfying 
$$p(A,B;C,X) = (A',B';C',P(X))\quad (=(P(A),P(B);P(C),P(X)).$$
The $P$ thus defined is injective: if $P(X)=P(Y)$, we have $p(A,B;C,X)= 
p(A,B;C,Y)$, and from injectivity of $p:G\to G'$ then follows that
$(A,B;C,X)= 
(A,B;C,Y)$. Proposition \ref{crrxx} for $\LL$ then gives that $X=Y$.
- Since $(P,p)$ thus preserves cross ratios of the form $(A,B;C,X)$, it 
follows from Proposition \ref{pex} that it also 
preserves all cross ratios of permutation instances thereof, i.e.\ 
cross ratios of 4-tuples, three of whose entries are $A,B$, and $C$.
But then also $(P,p)$ preserves all cross ratios of the form $(A,B;X,Y)$; 
for, by Proposition \ref{assxx}, 
$(A,B;X,Y)=(A,B;X,C).(A,B;C,Y)$, and since $p$ preserves the two 
cross ratios on the right, it preserves their composite, being a 
group homomorphism $G\to G'$. 
 From Proposition \ref{pex} now follows that $(P,p)$ preserves 
all cross ratios, two of whose entries are $A$ and $B$. 
But then also $(P,p)$ preserves all cross ratios of the form $(A,X;Y,Z)$, 
using
$ (A,X;Y,Z)=(A,X;Y,B).(A,X;B,Z)$, and again therefore all cross 
ratios one of whose entries is $A$. And then finally
by $(X,Y;Z,U)= (X,Y;Z,A).(X,Y;A,U)$, we see that $(P,p)$ preserves all 
cross ratios. 

We have  defined the value of $P$ on all objects of  $\LL$. 
To define the values of $P$ on the arrows of $\LL$, we put
$$P(\begin{diagram}A&\rTo^{C}& B\end{diagram})
:=\begin{diagram}P(A)&\rTo^{P(C)}& P(B)\end{diagram}$$ for arrows $A\to B$ with $A\neq 
B$; this then also ensures that $P$ preserves the projection structure. Arrows $A\to A$ may canonically be  identified with 
(abstract) scalars $\in G$, and the value on on such are 
given canonically by the assumed $p:G\to G'$. Henceforth, we don't 
distinguish notationally between $P$ and $p$. It remains to prove 
that $P$ preserves composition of arrows.
  
The $P$ constructed commutes with composition (multiplication) of 
scalars, by assumption on $p$. Also, it commutes with composition of 
non-scalars; for, assume that $(E:A\to B).(F:B\to C)= (G:A\to C)$, or 
equivalently, by Axiom 1, that
$$\begin{diagram}A&\rTo^{E}&B&\rTo^{F}&C&\rTo^{G}&A
\end{diagram}\quad  = \quad  1_{A}.$$ This means that the trirapport
$(A,B,C;E,F,G)$ equals $1$.  Now since $\LL$ has at least five points,
we may use Proposition \ref{bitrixx} to conclude that $P$ preserves 
tri-rapports, and so it follows that 
$$(P(A),P(B),P(C);P(E),P(F),P(G)) =1,$$ from which we again conclude 
from 
Axiom 1, and from the fact that $P$ is compatible the projection structure,
 that $$P(E:A\to B).P(F:B\to C)= P(G:A\to C),$$
as claimed.
It remains to be argued that $P$ commutes compositions of the form 
$\mu.(C:A\to B)$, with $\mu$ a scalar (and similarly with scalars 
multiplied on the right). 
 But we may  pick a point 
$D$ distinct from $A$, $B$,  $E$, and $F$ such that $\mu 
=(A,D;E,F)$; then  the composite $\mu.(C:A\to B)$ is
$$\begin{diagram}A&\rTo^{E}&D&\rTo^{F}&A&\rTo^{C}&B
\end{diagram},$$
with $\mu$ being the composite of the two leftmost arrows. But 
rebracketing, one gets a composite of two non-scalars $A\to D$ and 
$D\to B$, and the $D\to B$ in question  in turn is a composite of two 
non-scalars; so the whole three-fold product is preserved.

\section{Coordinatization}\label{coox}
 Given a field $k$, we get a $\Phi$-group $G$ 
by taking $G$ to be the group of non-zero elements, and taking $\Phi 
(\mu ) = 1-\mu$ for $\mu \neq 1$. This $\Phi$-group also has a 
specified element of order $\leq 2$, 
namely $-1$. A morphism of fields $k\to k'$ (is injective and) 
induces a morphism $G\to G'$ of the corresponding $\Phi$-groups, and 
it  also preserves $-1$. It is in fact a functor from the category of 
fields to the category of $\Phi$-groups with specified $-1$.
\begin{prop}This functor is full and faithful.
\end{prop}  
\pf  It is clear that the functor is faithful. To see that 
it is full, let $p:G\to G'$ be a  morphism of $\Phi$-groups with a 
specified 
$-1$. If $G$ and $G'$ come from fields $k$ and $k'$, respectively, 
the map $p\cup\{0\}: G\cup\{0\}\to G'\cup\{0\}$ is in fact a ring homomorphism: 
it clearly preserves multiplication, and it preserves addition, since 
addition can be reconstructed from the multiplication, $\Phi$ and 
$-1$, by the formula
\begin{equation}\label{xplusxx}\lambda + \mu= \lambda .\Phi ((-1).\lambda 
^{-1}. \mu)
\end{equation}
(and $0+\lambda =\lambda$, $0.\lambda = 0$).

\medskip

In other words, it makes sense to ask whether for a $\Phi$-group $G$ 
with specified $-1$, 
the algebraic structure $G\cup\{0\}$ (with the described $+$ and $.$ 
etc.) is a {\em field} or not.

To give the following Coordinatization Theorem a more succinct 
formulation, we pose the following definition.
Let $\LL$ be a groupoid with projection structure, satisfying Axioms 
1-4; let its $\Phi$-group of scalars with $-1$ be $G$.
Then $\LL$ is called a {\em 
projective line groupoid} if $G\cup\{0\}$ (with the addition $+$ given by 
(\ref{xplusxx}) etc.) is a field $k$ (called the {\em scalar field of} $\LL$). 

Then we have
\begin{thm} Every projective line groupoid $\LL$ with at least 
five\footnote{The result also holds in the case where $\LL$ has 
precisely 
three points, or if it has precisely four points and $-1\neq1$, see 
the Remark at the end of Section \ref{threefourx}. }
points  is isomorphic (as a 
groupoid with projection structure) to the groupoid $\LL (k^{2})$,
where $k$ is the field $G\cup \{0\}$.
\end{thm}
\pf  We know already from Sections \ref{sect1x} and 2 that $\LL (k^{2})$ is a 
projective line groupoid with scalar field $k$. Pick three distinct 
points $A,B,C$ of $\LL$, and three distinct points $A',B',C'$ of $\LL 
(k^{2})$. With $p$ as the the identity map $k\to k$, we can apply 
Theorem \ref{main3x} 
(Fundamental Theorem on three-transitivity) to obtain an morphism of 
projective line groupoids $P: \LL \to \LL (k^{2})$. It is surjective 
on objects by the last assertion in the  Theorem \ref{main3x}
 
\section{The twelve  scalars}
It is well known that permuting the four entries in a cross ratio 
expression $(A,B;C,D)$ with value $\mu$ gives six classical scalars:
$\mu, \mu^{-1}$, $1-\mu $, $(1-\mu )^{-1}$, $1 - \mu^{-1}$, and 
$(1-\mu^{-1})^{-1}$, see e.g.\ [Struik, 1953] p.\ 8. 
These relationships hold  in the axiomatic setting, without the 
further assumptions needed for the coordinatization, in other words, 
they hold if $1-x$ replaced  by $\Phi (x)$ everywhere (recall the involution 
$\Phi$ defined in terms of middle four interchange of cross ratio 
expressions). 

With the availability of the scalar $-1$, each 
of six classical scalars   derived from $\mu = (A,B;C,D)$,  
 by permutation of $A,B,C,D$, may be multiplied by $-1$, 
so that we get six further 
scalars. These, however, are not in general expressible by  cross 
ratios (= bi-rapports) built from the four given points $A,B,C,D$, but they are expressible as
tri-rapports bulit from them, using tri-rapport expressions for the six classical 
scalars. We refer to [Kock 2010b] for the full  list of the twelve 
tri-rapport expressions, but we include for 
convenience the basic tool for compiling this list:
\begin{prop}Let $A,B,C,D$ be mutually distinct. Then $$(A,B;C,D)= 
(A,C,D;B,A,B),$$
and $$(A,B,C;D,A,D) = - (A,B,D;D,A,C).$$
\end{prop}
\pf  In the diagram expressing Axiom 2, insert a North-East arrow 
$B:D\to A$ from  $D$ in the lower 
row  to the upper right hand $A$. The diagram then decomposes into a 
triangle and a pentagon. 
The triangle commutes by Axiom 1, and hence the 
pentagon commutes as well, and it expresses the relationship 
between the bi-rapport and the tri-rapport claimed in the first 
equation. The second equation was proved in Proposition 
\ref{signchangex}.

\section{An open end} Of particular geometric interest is projective 
geometry over the geometric line, often identified with the field 
${\mathbb R}$ of real numbers. However, it can be argued that other 
models $R$ of the geometric line exist, making differential geometry 
over manifolds based on $R$ more combinatorial or synthetic in 
nature, see e.g.\ [Kock, 1981]. But such an $R$ is not a field in the sense 
that the dichomtomy: ``for any $x\in R$, $x$ is either 0 or 
invertible''; rather $R$ is a 
local ring (inside a topos, in fact).

From this point of view, the present note is unsatisfactory: it 
depends heavily on the dichotomy -- here in the disguise that the 
arrows of the groupoid $\LL$ are {\em either} endomaps (scalars), 
{\em or} of the 
form $C:A\to B$ with $A,B$ and $C$ 
distinct objects.

Grassmanninan manifolds, say projective spaces like $P(R^{2})$, do 
exist and behave 
nicely when based on a local ring $R$, even in a topos , cf.\ e.g.\ 
[Kock-Reyes, 1977] or [Kock, 2010c],  A.5. However, our present axiomatics is too strong to cover 
such manifolds.

 An axiomatics that applies to this case is missing. 
For one thing,  the $\LL$ in question should not be a groupoid; even 
over a field $k$, some of the formulations given in the present note 
could be made more homogeneous by including {\em all} linear maps 
between 1-dimensional linear subspaces of $k^{2}$, not just the 
invertible ones, in other words, by including the zero-maps.
 In this case, $\LL$ satisfies the dichotomy:
``for any two objects $A$ and $B$, there is {\em exactly one} 
non-invertible map $A\to B$, and {\em at least one} invertible map $A\to 
B$''.

\bigskip

\noindent \verb+kock@imf.au.dk +


\begin{thebibliography}{99}

\bibitem{}  R. Baer (1952),  Linear Algebra and Projective Geometry, Academic Press 
1952 (Dover republication 2005).

\bibitem {}J.-L. Cathelineau (1995), Birapport et groupoides, L'Enseignement 
Math\'{e}matique 41 (1995), 257-280.

\bibitem{} H.S.M. Coxeter (1942), Non-Euclidean Geometry, University of Toronto 
Press 1942 (Fifth Ed.\ 1965).

\bibitem{} Y.\ Diers and J.\ Leroy (1994), Cat\'{e}gorie des point d'un espace 
projectif, Cahiers de Topologie et 
G\'{e}om\'{e}trie Diff\'{e}rentielle Cat\'{e}goriques 35 (1994), 2-28.

\bibitem{} A.\ Kock (1974), The category aspect of projective space, Aarhus 
Preprint Series 1974/75 No.\ 7.

\bibitem{} A.\ Kock,  Synthetic Differential Geometry, Cambridge 
University Press 1981; 2nd ed.\ 2006.

\bibitem{} A.\ Kock (2010a), Abstract projective lines, Cahiers de Topologie et 
G\'{e}om\'{e}trie Diff\'{e}rentielle Cat\'{e}goriques 51 (2010), 
224-240.


\bibitem{}  A.\ Kock (2010b), Geometric algebra of projective lines, 
arXiv:1003.2095, 2010.

\bibitem{}  A.\ Kock (2010c), Synthetic Geometry of Manifolds, Cambridge 
Tracts in Mathematics 180 (2010).

\bibitem{} A.\ Kock and G.E.\ Reyes (1977), Manifolds in formal differenial 
geometry, in ``Applications of Sheaves, Proceedings Durham 1977'', 
ed. M.\ Fourman et al., Springer Lecture Notes 753 (1979).


\bibitem{} D.\ Struik (1953), Analytic and Projective Geometry, 
Addison-Wesley 1953.

\end{thebibliography}
\end{document}